\documentclass[12pt]{article}
\usepackage{amssymb,amsmath,amsthm,amsfonts,enumerate, bbm, mathdots}

\usepackage{latexsym}
\usepackage{amscd}
\usepackage{stmaryrd}
\usepackage{color}
\usepackage[all]{xypic}
\usepackage{epsfig}
\usepackage{graphics}
\usepackage{ifthen}
\usepackage{varioref}
\usepackage{rotating}
\usepackage{extarrows}
\usepackage{cite}
\usepackage{mathrsfs}

\numberwithin{equation}{section}

\theoremstyle{plain}
\newtheorem{thm}{Theorem}[section]
\newtheorem{cor}[thm]{Corollary}
\newtheorem{lem}[thm]{Lemma}

\newtheorem{rem}[thm]{Remark}

\newcommand{\Res}{{\rm Res}}

\newcommand{\si}{\sigma}

\newcommand{\bfs}{{\boldsymbol{\sl{s}}}}

\definecolor{darkgreen}{rgb}{0.0625,0.64,0.0625}
\usepackage[
pdfstartview=FitH,
CJKbookmarks=true,
bookmarksnumbered=true,
bookmarksopen=true,
colorlinks,
pdfborder=001,
linkcolor=blue,
anchorcolor=green,
citecolor=red
]{hyperref}

\renewcommand{\Re}{\operatorname{Re}}
\newfont{\scyr}{wncyr10 scaled 550}

\def\proof{\noindent {\bf Proof.\;}}

\def\N{\mathbb{N}}
\def\Z{\mathbb{Z}}

\def\t{\widetilde{t}}

\def\ze{\zeta}

\allowdisplaybreaks

\begin{document}
	
\title{Parametric Euler $T$-sums of odd harmonic numbers}
\author{
{Ce Xu${}^{\text{a},}$\thanks{Email: cexu2020@ahnu.edu.cn.}\ ~and~ Lu Yan${}^{\text{b},}$\thanks{Email: 1910737@tongji.edu.cn, corresponding author.}}\\[1mm]
\small a. School of Mathematics and Statistics, \\ \small  Anhui Normal University, Wuhu 241002, China\\
\small b. School of Mathematical Sciences, \\ \small Tongji University, Shanghai 200092, China}
	
\date{ }

\maketitle	
	
\begin{abstract}
In this paper, we define a parametric variant of generalized Euler sums and call them the (alternating) parametric Euler $T$-sums. By using the contour integration method and residue theorem, we establish several explicit formulae for the linear parametric Euler $T$-sums. Furthermore, by applying the results, we obtain explicit formulae for the Hoffman's (alternating) double $t$-values and Kaneko-Tsumura's (alternating) double $T$-values.
\end{abstract}
	
	
	
	
\section{Introduction}\label{Sec:Intro}

We begin with some basic notations. Let $\mathbb{C},\mathbb{Z},\mathbb{N}$ and $\mathbb{N}^-$ be the sets of complex numbers, integers, positive integers and negative integers, respectively. We also denote by $\mathbb{N}_0$ the set of non-negative integers and by $\mathbb{N}^-_0$ the set of non-positive integers.

For $n,p\in\mathbb{N}$, the $n$-th generalized harmonic number of order $p$, denoted by $H_n^{(p)}$, is defined by
\begin{align}
H_n^{(p)}:=\sum\limits_{k=1}^{n}\frac{1}{k^p}.
\end{align}
Then $H_n^{(1)}\equiv H_n$ is the $n$-th harmonic number. For any $p\in\mathbb{N}$, we set $H_0^{(p)}:=0$. If $p>1$, the generalized harmonic number $H_n^{(p)}$ converges to the Riemann zeta value $\zeta(p)$ when $n$ tends to infinity, that is
$$\lim\limits_{n\rightarrow\infty}H_n^{(p)}=\zeta(p).$$
Here when $\Re{(z)}>1$, the Riemann zeta function $\zeta(z)$ is defined by
$$\zeta(z):=\sum_{n=1}^\infty\frac{1}{n^z}.$$

The linear Euler sum is defined by
\begin{align}
S_{p,q}:=\sum_{n=1}^\infty\frac{H_n^{(p)}}{n^q},
\end{align}
where $p,q\in\mathbb{N}$ with $q\geq2$.
It is known from  \cite[\citen{B-B-G}, \citen{Flajolet-Salvy}]{BBG} that the linear sums can be evaluated in terms of the Riemann zeta values in the following cases: $p=1$; $p=q$; $p+q$ odd and $p+q=6$ with $q\geq2$. In 1998, Flajolet and Salvy introduced the generalized Euler sum $S_{\boldsymbol{p},q}$ in \cite{Flajolet-Salvy}, which is defined as
\begin{align}\label{Defn-Cl-ES}
S_{\boldsymbol{p},q}:=\sum_{n=1}^\infty\frac{H_n^{(p_1)}H_n^{(p_2)}\cdots H_n^{(p_r)}}{n^q},
\end{align}
where $\boldsymbol{p}=(p_1, p_2,\ldots,p_r)\in\mathbb{N}^r$ and $q\in\mathbb{N}$ with $p_1\leq p_2\leq\cdots\leq p_r$ and $q\geq2$. In \cite{Flajolet-Salvy}, Flajolet and Salvy used the method of contour integration to evaluate the linear and generalized Euler sums. In \cite{Xu}, the first author of the paper used the same way to obtain some explicit evaluations of the parametric Euler sums.

The generalized Euler sums are in contrast to the multiple zeta values (abbr. MZVs) defined as
\begin{align}\label{Defn-MZV}
\zeta(\bfs)\equiv \zeta(s_1,\ldots,s_k):=\sum\limits_{n_1>\cdots>n_k>0 } \frac{1}{n_1^{s_1}\cdots n_k^{s_k}},
\end{align}
where $\bfs=(s_1,s_2,\ldots,s_k)\in \N^k$ with $s_1>1$. The systematic study of MZVs began in the early 1990s with the works of Hoffman \cite{H1992} and Zagier \cite{DZ1994}. Due to surprising applications in many branches of mathematics and theoretical physics, MZVs have attracted a lot of attention and interest in the past three decades (for example, see the books by Srivastava-Choi \cite{SC2012} and Zhao \cite{Zhao2016}).  As Euler \cite{Euler} discovered by a process of extrapolation, the double zeta values $\zeta(s,t)$ can be evaluated in terms of the Riemann zeta values when $s+t$ is odd. Later, Borweins and Girgensohn\cite{B-B-G} obtained explicit formulae for $\zeta(s,t)$. In the study of some parametric Euler sums, Borweins and Bradley deduced an explicit formula of the double zeta values $\zeta(2j, 2m+1)$ by the method of power series expansion and comparing coefficients in \cite{BBB}. For more results about generalized Euler sums, please see \cite[\citen{BBB}, \citen{C2011}-\citen{C2013}, \citen{Mezo}, \citen{S2011}-\citen{SS2011}, \citen{W-L}, \citen{Xu2}, \citen{Xu-Wang}]{AC} and the references therein.

Recently, an odd variant of the Euler sums was introduced in \cite{Xu-Wang2}. For $\boldsymbol{p}=(p_1, p_2,\ldots,p_r)\in\mathbb{N}^r$ and $q\in\mathbb{N}$ with $p_1\leq p_2\leq\cdots\leq p_r$ and $q\geq2$, the authors of \cite{Xu-Wang2} defined
\begin{align}
T_{\boldsymbol{p},q}:=\sum_{n=1}^\infty\frac{h_{n-1}^{(p_1)}h_{n-1}^{(p_2)}\cdots h_{n-1}^{(p_r)}}{\left(n-\frac{1}{2}\right)^q},
\end{align}
which is called an Euler $T$-sum. Here for $n,p\in\mathbb{N}$, $h_n^{(p)}$ is the $n$-th odd harmonic number of order $p$ defined by
$$h_n^{(p)}:=\sum\limits_{k=1}^{n}\frac{1}{(k-\frac{1}{2})^p},$$
and $h_0^{(p)}=0$. They further considered the alternating Euler $T$-sums
\begin{align}
T_{\boldsymbol{p},\overline{q}}:=\sum_{n=1}^\infty(-1)^n\frac{h_{n-1}^{(p_1)}h_{n-1}^{(p_2)}\cdots h_{n-1}^{(p_r)}}{\left(n-\frac{1}{2}\right)^q},
\end{align}
where $q$ can be any positive integer. In \cite{Xu-Wang2}, the authors obtained the explicit formulae of linear $T$-sums $T_{p,q}$ with $p+q$ odd, quadratic $T$-sums $T_{p_1p_2,q}$ with $p_1+p_2+q$ even and cubic $T$-sums $T_{1^3,q}$ with $q$ even. They also establish explicit formulae for the alternating linear $T$-sums $T_{p,\bar{q}}$ and alternating quadratic $T$-sums $T_{1^2,\bar{q}}$. Furthermore, they concluded that the triple $t$-values are reducible to the Riemann zeta values, double zeta values and double $t$-values, and the triple $T$-values $T(s_1,s_2,s_3)$ with $s_1+s_2+s_3$ even can be expressed in terms of single and double $T$-values.
Here for $s_1,s_2,\ldots,s_k\in\N$ with $s_1>1$, the multiple $t$-value $t(s_1,s_2,\ldots,s_k)$ introduced by Hoffman in \cite{Hoffman}, and the multiple $T$-value $T(s_1,s_2,\ldots,s_k)$ introduced by Kaneko and Tsumura in \cite{KTA2018,KTA2019} are defined as
\begin{align}\label{MtVs}
	t(s_1,s_2,\ldots,s_k):=\sum\limits_{n_1>n_2>\cdots>n_k\geq1}{\frac{1}{(2n_1-1)^{s_1}(2n_2-1)^{s_2}\cdots(2n_k-1)^{s_k}}},
\end{align}
and
\begin{align}
	T(s_1,s_2,\ldots,s_k)
	&:=2^k\sum_{{n_1>n_2>\cdots>n_k\geq1\atop n_i\equiv k-i+1\ {\rm mod}\ {2}}}
	\frac{1}{n_1^{s_1}n_2^{s_2}\cdots n_k^{s_k}}\nonumber\\
	&=2^k\sum_{n_1>n_2>\cdots>n_k\geq 1}
	\frac{1}{(2n_1-k)^{s_1}(2n_2-k+1)^{s_2}\cdots(2n_k-1)^{s_k}}\,,\label{MTV}
\end{align}
respectively. For $(s_1,s_2,\ldots,s_k)\in\mathbb{N}^k$ and $(\si_1,\si_2,\ldots,\si_k)\in\{\pm1\}^k$ with $(s_1,\si_1)\neq (1,1)$, we define the alternating multiple $t$-values and alternating multiple $T$-values by
\begin{align}\label{AMtVs}
	t(s_1,s_2,\ldots,s_k;\si_1,\si_2,\ldots,\si_k)
	&:=\sum_{n_1>n_2>\cdots>n_k\geq 1}\frac{\si_1^{n_1}\si_2^{n_2}\cdots\si_k^{n_k}}
	{(2n_1-1)^{s_1}(2n_2-1)^{s_2}\cdots(2n_k-1)^{s_k}}\,,
\end{align}
and
\begin{align}\label{AMTVs}
	&T(s_1,s_2,\ldots,s_k;\si_1,\si_2,\ldots,\si_k)\nonumber\\
	&\quad:=2^k\sum_{n_1>n_2>\cdots>n_k\geq 1}
	\frac{\si_1^{n_1}\si_2^{n_2}\cdots\si_k^{n_k}}
	{(2n_1-k)^{s_1}(2n_2-k+1)^{s_2}\cdots(2n_k-1)^{s_k}}\,,
\end{align}
respectively. Throughout this paper, we may denote the alternating multiple $t$-values and alternating multiple $T$-values more concisely by a rule:
whenever $\sigma_j=-1$,  we place a bar over the corresponding component $s_j$. For example, we have
\begin{equation*}
	t(\overline{s_1}, s_2)=t(s_1,s_2;-1,1)\quad \text{and}\quad T(\overline{s_1},s_2,\overline{s_3})=T(s_1,s_2,s_3;-1,1,-1).
\end{equation*}

In this paper, we consider the following (alternating) parametric Euler sums involving odd harmonic numbers:
\begin{align}\label{studyobject}
	T_{\boldsymbol{p},\boldsymbol{q}}^\si(a_1,a_2,\ldots,a_k):=\sum_{n=1}^\infty\frac{h_n^{(p_1)}h_n^{(p_2)} \cdots h_n^{(p_r)}\si^n}{\left(n+a_1-\frac{1}{2}\right)^{q_1}\left(n+a_2-\frac{1}{2}\right)^{q_2}\cdots\left(n+a_k-\frac{1}{2}\right)^{q_k}},
 \end{align}
and call them the (alternating) parametric Euler $T$-sums. Here $\si\in\{\pm1\}$, $a_1-\frac{1}{2},a_2-\frac{1}{2},\ldots,a_k-\frac{1}{2}\notin\mathbb{N}^-$, $\boldsymbol{p}=(p_1, p_2,\ldots,p_r)\in\mathbb{N}^r$ with $p_1\leq p_2\leq\cdots\leq p_r$, $\boldsymbol{q}=(q_1, q_2,\ldots,q_k)\in\mathbb{N}_0^k$ with $q_1+q_2+\cdots+q_k\geq2$ for $\si=1$ and $q_1+q_2+\cdots+q_k\geq1$ for $\si=-1$. We will give some formulae for the (alternating) parametric linear Euler $T$-sums by applying the contour integral representation and residue theorem. As applications, similarly as in \cite{BBB}, we show that the double $t$-values $t(s_1,s_2)$ can be evaluated in terms of the Riemann zeta values and single $t$-values when $s_1+s_2$ is odd and the alternating double $t$-values $t(\overline{s_1},s_2)$ can be evaluated in terms of the alternating Riemann zeta values, alternating single $t$-values and single $t$-values when $s_1+s_2$ is even. Moreover, we get similar conclusions for  the double $T$-values $T(s_1,s_2)$ with $s_1+s_2$ odd and the alternating double $T$-values $T(\overline{s_1},s_2)$ with $s_1+s_2$ odd.

The plan of the paper is as follows. In Section 2 we introduce the residue lemma and give some asymptotic formulae for the parametric digamma functions and the trigonometric functions. In Section 3, we use the approach of contour integral to study the (alternating) parametric linear Euler $T$-sums. In the last section, by applying some explicit evaluations obtained in Section 3, we prove explicit formulae for double $t$-values, double $T$-values, alternating double $t$-values and alternating double $T$-values.

	
\section{Several Lemmas}\label{Sec:lemmas}
	
A complex function $\xi(z)$ is called a kernel function if
\begin{itemize}
  \item [(i)] $\xi(z)$ is meromorphic in the whole complex plane,
  \item [(ii)] $\xi(z)$ satisfies $\xi (z)=o(z)$ over an infinite collection of circles $|z|=\rho_k$ with $\rho_k\to \infty $.
\end{itemize}
 Applying these two conditions, Flajolet and Salvy \cite{Flajolet-Salvy} discovered the following residue lemma.

\begin{lem}[\cite{Flajolet-Salvy}]\label{residue lemma}
Let $\xi(z)$ be a kernel function and let $r(z)$ be a rational function which is $O(z^{-2})$ at infinity. Then
\begin{align*}
\sum_{\alpha\in O} \Res(r(z)\xi(z),\alpha)+ \sum_{\beta\in S} \Res(r(z)\xi(z),\beta) = 0,
\end{align*}
where $S$ is the set of poles of $r(z)$ and $O$ is the set of poles of $\xi(z)$ that are not poles of $r(z)$. Here $\Res(f(z),\lambda)$ denotes the residue of $f(z)$ at $z=\lambda$.
\end{lem}

\begin{lem}[{\cite[Lemma 2.2]{Xu-Wang2}}]\label{lemma2}
If a meromorphic function $f(z)$ has a pole of order $m$ at $z=\lambda$, then
$$\Res\left(f(z), \lambda\right)
 =\lim\limits_{z\rightarrow\lambda}\frac{1}{(m-1)!}\frac{d^{m-1}}{ds^{m-1}}\left[(z-\lambda)^mf(z)\right]
  =\lim\limits_{z\rightarrow\lambda}\frac{1}{m!}\frac{d^m}{ds^m}\left[(z-\lambda)^{m+1}f(z)\right].$$
\end{lem}

The parametric digamma function, denoted by $\boldsymbol{\Psi}(-z;a)$ in \cite{Xu2}, is defined by
$$\boldsymbol{\Psi}(-z;a)+\gamma
	=\frac{1}{z-a}+\sum\limits_{k=1}^{\infty}\left(\frac{1}{k+a}-\frac{1}{k+a-z}\right),$$
where  $a\notin\mathbb{N}^-$ and $\gamma$ denotes the Euler-Mascheroni constant. The function $\boldsymbol{\Psi}(-z;a)$ is meromorphic in the entire complex plane with a simple pole at $z=n+a$ for each negative integer $n$. For simplicity, we let
$$\boldsymbol{\Psi}(-z)=\boldsymbol{\Psi}\left(-z;-\frac{1}{2}\right)+\gamma
=\frac{1}{z+\frac{1}{2}}+\sum\limits_{k=1}^{\infty}\left(\frac{1}{k-\frac{1}{2}}-\frac{1}{k-\frac{1}{2}-z}\right).$$
Obviously,  the $(p-1)$-th derivative of the function $\boldsymbol{\Psi}(-z)$ is
\begin{align}\label{digamma}
		\boldsymbol{\Psi}^{(p-1)}(-z)=(-1)^p(p-1)!\sum\limits_{k=0}^{\infty}\frac{1}{\left(k-\frac{1}{2}-z\right)^p}.
\end{align}

To introduce the following lemmas more concisely, we define
	\begin{align}\label{widetilde-t}
	     \widetilde{t}(s_1,s_2,\ldots,s_k)&:=2^{s_1+s_2+\cdots+s_k}t(s_1,s_2,\ldots,s_k)\notag\\
	     &=\sum\limits_{n_1>n_2>\cdots>n_k\geq1}{\frac{1}{\left(n_1-\frac{1}{2}\right)^{s_1}\left(n_2-\frac{1}{2}\right)^{s_2}\cdots\left(n_k-\frac{1}{2}\right)^{s_k}}},
	\end{align}
    and call them multiple $\widetilde{t}$-values. Similarly, the corresponding alternating multiple $\widetilde{t}$-values are defined by
    \begin{align}
    \t(s_1,s_2,\ldots,s_k;\si_1,\si_2,\ldots,\si_k)
    &:=\sum_{n_1>n_2>\cdots>n_k\geq 1}\frac{\si_1^{n_1}\si_2^{n_2}\cdots\si_k^{n_k}}
        {(n_1-1/2)^{s_1}(n_2-1/2)^{s_2}\cdots(n_k-1/2)^{s_k}}\nonumber\\
    &=2^{s_1+s_2+\cdots+s_k}t(s_1,s_2,\ldots,s_k;\si_1,\si_2,\ldots,\si_k).
    \end{align}
    Similarly, we put a bar on the top of $s_j$ if $\si_j=-1$. For example, we have
    \begin{align*}
    	\widetilde{t}(\overline{s})= \widetilde{t}(s;-1)
        \quad \text{and}\quad
         \widetilde{t}(\overline{s_1}, s_2)=\widetilde{t}(s_1,s_2;-1,1)
    \end{align*}

    From \cite[Theorems 2.1-2.3, Corollary 2.4]{Xu2} and \cite[(2.2)-(2.7)]{Xu-Wang2}, we list the following asymptotic formulae of the function  $\boldsymbol{\Psi}^{(p-1)}\left(\frac{1}{2}-z\right)$ at the integers and poles.
   	\begin{lem}\label{parametic-digamma}
    	For $n\in\N_0$ and $p\in\mathbb{N}$, we have
    \begin{align}
       &\begin{aligned}\label{psi n}
    	\frac{\boldsymbol{\Psi}^{(p-1)}\left(\frac{1}{2}-z\right)}{(p-1)!}\overset{z\rightarrow n}{=}
    	\frac{1}{(z-n)^p}+(-1)^p\sum\limits_{i=p}^{\infty}\dbinom{i-1}{p-1}\left((-1)^iH_n^{(i)}+\zeta(i)\right)(z-n)^{i-p},
        \end{aligned}\\
       &\begin{aligned}\label{psi n-1/2}
    	\frac{\boldsymbol{\Psi}^{(p-1)}\left(\frac{1}{2}-z\right)}{(p-1)!}\overset{z\rightarrow n-\frac{1}{2}}{=}
    	(-1)^p\sum\limits_{i=p}^{\infty}\dbinom{i-1}{p-1}\left((-1)^ih_n^{(i)}+\widetilde{t}(i)\right)\left(z-n+\frac{1}{2}\right)^{i-p},
        \end{aligned}\\
       &\begin{aligned}\label{psi 1/2-n}
    	\frac{\boldsymbol{\Psi}^{(p-1)}\left(\frac{1}{2}-z\right)}{(p-1)!}\overset{z\rightarrow \frac{1}{2}-n}{=}
    	(-1)^p\sum\limits_{i=p}^{\infty}\dbinom{i-1}{p-1}\left(\widetilde{t}(i)-h_{n-1}^{(i)}\right)\left(z+n-\frac{1}{2}\right)^{i-p},\quad
    	n\in\mathbb{N},
        \end{aligned}
        \end{align}
   where $\ze(1):=-2\log(2)$ and $\t(1):=0$.
    \end{lem}

    Clearly,
   $$\boldsymbol{\Psi}\left(\frac{1}{2}-z\right)=\psi(-z)-\psi(1/2),\quad \psi(1/2)=-\gamma-2\log(2).$$
   Here $\psi (z)$ is the classical digamma function defined for $z\notin  \N^-_0$ by
    \[\psi \left(z \right) =  - \gamma  + \sum\limits_{n = 0}^\infty  {\left( {\frac{1}{{n + 1}} - \frac{1}{{n + z}}} \right)}.\]

According to \cite[\citen{Wang-Guo}, \citenum{Xu2}, \citenum{Xu-Wang2}]{AAR}, the following expressions of $\pi \tan(\pi z)$ and $\frac{\pi}{\cos(\pi z)}$ hold.

\begin{lem}
For any $n\in \Z$ and $p\in\N_0$, we have
    \begin{align}
	    &\begin{aligned}\label{tan n-1/2}
	    	\pi \tan(\pi z)\overset{z\rightarrow n-\frac{1}{2}}{=}
	    	-\frac{1}{z-n+\frac{1}{2}}+2\sum\limits_{i=1}^{\infty}\zeta(2i)\left(z-n+\frac{1}{2}\right)^{2i-1},
	    \end{aligned}\\
	    &\begin{aligned}\label{tan}
			\pi \tan(\pi z)=
			2\sum\limits_{i=1}^{\infty}\widetilde{t}(2i)z^{2i-1},\quad |z|<1,
		\end{aligned}\\
     	&\begin{aligned}\label{d-tan}
		    \lim\limits_{z\rightarrow n}\frac{d^p}{dz^p}(\pi \tan(\pi z))
		    =(1-(-1)^p)p!\widetilde{t}(p+1),
    	\end{aligned}\\
	    &\begin{aligned}\label{cos n-1/2}
		    \frac{\pi}{\cos(\pi z)}\overset{z\rightarrow n-\frac{1}{2}}{=}
		    (-1)^n\left\{\frac{1}{z-n+\frac{1}{2}}+2\sum\limits_{i=1}^{\infty}\overline{\zeta}(2i)\left(z-n+\frac{1}{2}\right)^{2i-1}\right\},
     	\end{aligned}\\
		&\begin{aligned}\label{cos}
			\frac{\pi}{\cos(\pi z)}=
			-2\sum\limits_{i=0}^{\infty}\widetilde{t}(\overline{2i+1})z^{2i},\quad |z|<1,
		\end{aligned}\\
		&\begin{aligned}\label{d-cos}
			\lim\limits_{z\rightarrow n}\frac{d^p}{dz^p}\frac{\pi}{\cos(\pi z)}
			=(-1)^{n-1}(1+(-1)^p)p!\widetilde{t}(\overline{p+1}),
		\end{aligned}
    \end{align}
		where $\overline{\zeta}(z)$ denotes the alternating Riemann zeta function which is defined for $\Re(z)\geq1$ by
		$$\overline{\zeta}(z)=\sum\limits_{n=1}^{\infty}\frac{(-1)^{n-1}}{n^z}.$$
\end{lem}


	\section{Parametric Euler $T$-sums}\label{Sec:proof}
	
	In this section, we apply the method of contour integration to study the (alternating) parametric linear Euler $T$-sums	
	$$T_{p, 1^2}(a, b)=\sum_{n=1}^{\infty}\frac{h_n^{(p)}}{\left(n+a-\frac{1}{2}\right)\left(n+b-\frac{1}{2}\right)}$$
	and
	$$T^{-1}_{p, 1^2}(a, b)=\sum_{n=1}^{\infty}(-1)^n\frac{h_n^{(p)}}{\left(n+a-\frac{1}{2}\right)\left(n+b-\frac{1}{2}\right)}.$$
	\begin{thm}\label{thm3.1}
		For any $p\in\mathbb{N}, a,b\in\mathbb{C}$ with $a\neq b, a,b\notin\mathbb{N}_0^{-}$ and $a+\frac{1}{2}, b+\frac{1}{2}\notin\mathbb{Z}$, we have
		\begin{align}\label{thm3.1-1}
			&\sum_{n=1}^{\infty}\frac{h_n^{(p)}}{\left(n+a-\frac{1}{2}\right)\left(n+b-\frac{1}{2}\right)}-(-1)^p\sum_{n=1}^{\infty}\frac{h_{n-1}^{(p)}}{\left(n-a-\frac{1}{2}\right)\left(n-b-\frac{1}{2}\right)}\notag\\
			=&2\frac{(-1)^p}{b-a}\sum\limits_{k=1}^{\left[\frac{p}{2}\right]}\widetilde{t}(2k)\left\{\zeta(p-2k+1; a)-\zeta(p-2k+1; b)\right\}\notag\\
			&+\frac{(-1)^p}{b-a}\left\{\pi \tan(\pi b)\left(\zeta(p;b)-\widetilde{t}(p)\right)-\pi \tan(\pi a)\left(\zeta(p;a)-\widetilde{t}(p)\right)\right\},
		\end{align}
    where $[x]$ is the integer part of $x$, and the Hurwitz zeta function $\zeta\left(z;a\right)$ is defined for $\Re\left(z\right)>1$ and $a\notin\mathbb{N}_0^{-}$ by
    $$\zeta\left(z;a\right):=\sum\limits_{n=0}^{\infty}{\frac{1}{\left(n+a\right)^z}}.$$
    We also set $\ze(1;a):=-\boldsymbol{\Psi}(1/2+a)=\psi(1/2)-\psi(a)$.
	\end{thm}
	\proof
	Let
	$$F(z)=\frac{\pi \tan(\pi z) \boldsymbol{\Psi}^{(p-1)}\left(\frac{1}{2}-z\right)}{(z+a)(z+b)(p-1)!}.$$
	This function has poles at $z=-a,-b,n(n\in\mathbb{N}_0)$ and $\pm\left(n-\frac{1}{2}\right) (n\in\mathbb{N})$.
	Obviously, $\frac{\pi \tan(\pi z) \boldsymbol{\Psi}^{(p-1)}\left(\frac{1}{2}-z\right)}{(p-1)!}$ is a kernel function and $\oint_{(\infty)} F(z)\mathrm{d}z=0$, where $\oint_{\left( \infty  \right)}$ denotes the integral along large circles, that is, the limit of integrals $\oint_{\left| z \right| = R}$ as $R\to \infty$. Hence, using Lemma \ref{residue lemma}, we obtain
	\begin{align}\label{thm3.1-2}
		&\sum\limits_{n=0}^{\infty}\Res(F(z), n)+\Res(F(z), -a)+ \Res(F(z), -b)\notag\\
		&+\sum\limits_{n=1}^{\infty}\left\{\Res\left(F(z), n-\frac{1}{2}\right)+\Res\left(F(z), \frac{1}{2}-n\right)\right\}=0.
	\end{align}
	Applying \eqref{psi n} and the fact that $z=n$ is a simple zero of $\tan(\pi z)$, we note that the pole of $F(z)$ at $z=n$ is of order $p-1$ for each nonnegative integer $n$. From Eqs. \eqref{psi n}, \eqref{d-tan} and Lemma \ref{lemma2}, the residue is
	\begin{align*}
		\Res(F(z), n)
		=&\frac{1}{(p-1)!}\lim\limits_{z\rightarrow n}\frac{d^{p-1}}{dz^{p-1}}\left\{(z-n)^pF(s)\right\}\\
		=&\frac{1}{(p-1)!}\lim\limits_{z\rightarrow n}\frac{d^{p-1}}{dz^{p-1}}\frac{\pi \tan(\pi z)}{(z+a)(z+b)}\\
		=&2\frac{(-1)^p}{b-a}\sum\limits_{k=1}^{\left[\frac{p}{2}\right]}\widetilde{t}(2k)\left(\frac{1}{(n+a)^{p-2k+1}}-\frac{1}{(n+b)^{p-2k+1}}\right).
	\end{align*}
	For a positive integer $n$, the poles of $F(z)$ at $z=\pm\left(n-\frac{1}{2}\right), -a, -b$ are simple and by using Eqs. \eqref{psi n-1/2}, \eqref{psi 1/2-n}, \eqref{tan n-1/2}, we deduce that
    \begin{align*}
	&\begin{aligned}
		\Res\left(F(z), n-\frac{1}{2}\right)
		=&-\lim\limits_{z\rightarrow n-\frac{1}{2}}\frac{\boldsymbol{\Psi}^{(p-1)}\left(\frac{1}{2}-z\right)}{(z+a)(z+b)(p-1)!}\\
		=&-\frac{(-1)^p\widetilde{t}(p)+h_n^{(p)}}{\left(n+a-\frac{1}{2}\right)\left(n+b-\frac{1}{2}\right)},
	\end{aligned}\\
	&\begin{aligned}
		\Res\left(F(z), \frac{1}{2}-n\right)
		=&-\lim\limits_{z\rightarrow \frac{1}{2}-n}\frac{\boldsymbol{\Psi}^{(p-1)}\left(\frac{1}{2}-z\right)}{(z+a)(z+b)(p-1)!}\\
		=&-(-1)^p\frac{\widetilde{t}(p)-h_{n-1}^{(p)}}{\left(n-a-\frac{1}{2}\right)\left(n-b-\frac{1}{2}\right)},
	\end{aligned}\\
	&\begin{aligned}
		\Res(F(z), -a)=-\frac{\pi \tan(\pi a) \boldsymbol{\Psi}^{(p-1)}\left(\frac{1}{2}+a\right)}{(b-a)(p-1)!},
	\end{aligned}\\
	&\begin{aligned}
	   \Res(F(z), -b)=\frac{\pi \tan(\pi b) \boldsymbol{\Psi}^{(p-1)}\left(\frac{1}{2}+b\right)}{(b-a)(p-1)!}.
    \end{aligned}
    \end{align*}
Then noting \eqref{digamma} and the fact that
	\begin{align*}
		&\frac{\pi}{b-a}\left(\tan(\pi b)-\tan(\pi a)\right)\\
		=&\sum\limits_{n=0}^{\infty}\frac{1}{\left(n+a+\frac{1}{2}\right)\left(n+b+\frac{1}{2}\right)}
		+\sum\limits_{n=1}^{\infty}\frac{1}{\left(n-a-\frac{1}{2}\right)\left(n-b-\frac{1}{2}\right)},
	\end{align*}
   we obtain \eqref{thm3.1-1} by \eqref{thm3.1-2} and a direct calculation.  So this completes the proof.
    \qed

	Taking $b=-a, p=2m+1$ with $m\in\mathbb{N}_0$ in \eqref{thm3.1-1}, we get the following corollary.
	\begin{cor}\label{cor1}
		For $m\in\mathbb{N}_0$ and $a\in\mathbb{C}$ with $a,a+\frac{1}{2}\notin\mathbb{Z}$, we have
		\begin{align}\label{cor11}
			&\sum_{n=1}^{\infty}\frac{h_n^{(2m+1)}}{\left(n+a-\frac{1}{2}\right)\left(n-a-\frac{1}{2}\right)}
            =\sum_{n=1}^{\infty}\frac{h_n^{(2m+1)}}{\left(n-\frac1{2}\right)^2-a^2}\notag\\
			=&\frac{1}{2}\sum_{n=1}^{\infty}\frac{1}{\left(n-\frac{1}{2}\right)^{2m+1}\left(\left(n-\frac1{2}\right)^2-a^2\right)}\notag\\
			&+\frac{1}{2a}\sum_{k=1}^m\widetilde{t}(2k)\left(\zeta(2m-2k+2; a)-\zeta(2m-2k+2; -a)\right)\notag\\
			&+\frac{1}{4a}\pi \tan(\pi a)\left(2\widetilde{t}(2m+1)-\zeta(2m+1; a)-\zeta(2m+1; -a)\right).
		\end{align}
	\end{cor}
    Putting $b=1-a, p=2m+1$ with $m\in\mathbb{N}_0$ in \eqref{thm3.1-1}, we obtain the following corollary.
    \begin{cor}\label{cor2}
	    For $m\in\mathbb{N}_0$ and $a\in\mathbb{C}$ with $a,a+\frac{1}{2}\notin\mathbb{Z}$, the following equality holds:
	    \begin{align}\label{cor22}
	    	&\sum_{n=1}^{\infty}\frac{h_n^{(2m+1)}}{\left(n+a-\frac{1}{2}\right)\left(n-a+\frac{1}{2}\right)}
            =\sum_{n=1}^{\infty}\frac{h_n^{(2m+1)}}{n^2-\left(a-\frac1{2}\right)^2}\notag\\
	    	=&\frac{1}{2a-1}\sum_{k=1}^m\widetilde{t}(2k)\left(\zeta(2m-2k+2; a)-\zeta(2m-2k+2; 1-a)\right)\notag\\
	    	&+\frac{\pi \tan(\pi a)}{2(2a-1)}\left(2\widetilde{t}(2m+1)-\zeta(2m+1; a)-\zeta(2m+1; 1-a)\right).
	    \end{align}
    \end{cor}

    \begin{thm}\label{thm3.4}
    	For any $p\in\mathbb{N}, a,b\in\mathbb{C}$ with $a\neq b, a,b\notin\mathbb{N}_0^{-}$ and $a+\frac{1}{2}, b+\frac{1}{2}\notin\mathbb{Z}$, we have
    	\begin{align}\label{thm3.4-1}
    		&\sum_{n=1}^{\infty}\frac{(-1)^nh_n^{(p)}}{\left(n+a-\frac{1}{2}\right)\left(n+b-\frac{1}{2}\right)}+(-1)^p\sum_{n=1}^{\infty}\frac{(-1)^nh_{n-1}^{(p)}}{\left(n-a-\frac{1}{2}\right)\left(n-b-\frac{1}{2}\right)}\notag\\
    		=&2\frac{(-1)^p}{b-a}\sum\limits_{k=0}^{\left[\frac{p-1}{2}\right]}\widetilde{t}(\overline{2k+1})\left\{\overline{\zeta}(p-2k; b)-\overline{\zeta}(p-2k; a)\right\}\notag\\
    		&+\frac{(-1)^p}{b-a}\left\{\frac{\pi}{\cos(\pi b)}\left(\zeta(p;b)-\widetilde{t}(p)\right)-\frac{\pi}{\cos(\pi a)}\left(\zeta(p;a)-\widetilde{t}(p)\right)\right\},
    	\end{align}
        where the alternating Hurwitz zeta function $\overline{\zeta}(z;a)$ is defined for $\Re(z)\geq1$ and $a\notin\mathbb{N}_0^{-}$ by
        $$\overline{\zeta}(z;a):=\sum\limits_{n=0}^{\infty}{\frac{(-1)^n}{\left(n+a\right)^z} }.$$
    \end{thm}
    \proof
    Let
    $$G(z)=\frac{\pi \boldsymbol{\Psi}^{(p-1)}\left(\frac{1}{2}-z\right)}{\cos(\pi z)(z+a)(z+b)(p-1)!}.$$
    This function has poles at $z=-a,-b,n(n\in\mathbb{N}_0)$ and $\pm\left(n-\frac{1}{2}\right) (n\in\mathbb{N})$.
    Clearly, $\frac{\pi \boldsymbol{\Psi}^{(p-1)}\left(\frac{1}{2}-z\right)}{\cos(\pi z)(p-1)!}$ is a kernel function and $\oint_{(\infty)} G(z)\mathrm{d}z=0$. By using Lemma \ref{residue lemma}, we have
    \begin{align}\label{thm3.4-2}
    	&\sum\limits_{n=0}^{\infty}\Res(G(z), n)+\Res(G(z), -a)+ \Res(G(z), -b)\notag\\
    	&+\sum\limits_{n=1}^{\infty}\left\{\Res\left(G(z), n-\frac{1}{2}\right)+\Res\left(G(z), \frac{1}{2}-n\right)\right\}=0.
    \end{align}
    For a nonnegative integer $n$, the pole of $G(z)$ at $z=n$ is of order $p$. From Eqs. \eqref{psi n}, \eqref{d-cos} and Lemma \ref{lemma2}, the residue is
    \begin{align*}
    	\Res(G(z), n)
    	=&\frac{1}{(p-1)!}\lim\limits_{z\rightarrow n}\frac{d^{p-1}}{dz^{p-1}}\left\{(z-n)^pG(s)\right\}\\
    	=&\frac{1}{(p-1)!}\lim\limits_{z\rightarrow n}\frac{d^{p-1}}{dz^{p-1}}\frac{\pi}{\cos(\pi z)(z+a)(z+b)}\\
    	=&2\frac{(-1)^{p+n}}{b-a}\sum\limits_{k=0}^{\left[\frac{p-1}{2}\right]}\widetilde{t}(\overline{2k+1})\left(\frac{1}{(n+a)^{p-2k}}-\frac{1}{(n+b)^{p-2k}}\right).
    \end{align*}
  For a positive integer $n$, the poles of $G(z)$ at $z=\pm\left(n-\frac{1}{2}\right), -a, -b$ are simple and by using Eqs. \eqref{psi n-1/2}, \eqref{psi 1/2-n}, \eqref{cos n-1/2}, we deduce that
    \begin{align*}
    &\begin{aligned}
    	\Res\left(G(z), n-\frac{1}{2}\right)
    	=&(-1)^n\lim\limits_{z\rightarrow n-\frac{1}{2}}\frac{\boldsymbol{\Psi}^{(p-1)}\left(\frac{1}{2}-z\right)}{(z+a)(z+b)(p-1)!}\\
    	=&(-1)^n\frac{(-1)^p\widetilde{t}(p)+h_n^{(p)}}{\left(n+a-\frac{1}{2}\right)\left(n+b-\frac{1}{2}\right)},
    \end{aligned}\\
    &\begin{aligned}
    	\Res\left(G(z), \frac{1}{2}-n\right)
    	=&(-1)^{n+1}\lim\limits_{z\rightarrow \frac{1}{2}-n}\frac{\boldsymbol{\Psi}^{(p-1)}\left(\frac{1}{2}-z\right)}{(z+a)(z+b)(p-1)!}\\
    	=&(-1)^{n+p+1}\frac{\widetilde{t}(p)-h_{n-1}^{(p)}}{\left(n-a-\frac{1}{2}\right)\left(n-b-\frac{1}{2}\right)},
    \end{aligned}\\
    &\begin{aligned}
    	\Res(G(z), -a)=\frac{\pi \boldsymbol{\Psi}^{(p-1)}\left(\frac{1}{2}+a\right)}{\cos(\pi a)(b-a)(p-1)!},
    \end{aligned}\\
    &\begin{aligned}
    	\Res(G(z), -b)=\frac{\pi \boldsymbol{\Psi}^{(p-1)}\left(\frac{1}{2}+b\right)}{\cos(\pi b)(a-b)(p-1)!}.
    \end{aligned}
      \end{align*}
 Now noting \eqref{digamma} and the fact that
    \begin{align*}
    	&\frac{\pi}{b-a}\left(\frac{1}{\cos(\pi a)}-\frac{1}{\cos(\pi b)}\right)\\
    	=&\sum\limits_{n=0}^{\infty}\frac{(-1)^n}{\left(n+a+\frac{1}{2}\right)\left(n+b+\frac{1}{2}\right)}
    	+\sum\limits_{n=1}^{\infty}\frac{(-1)^n}{\left(n-a-\frac{1}{2}\right)\left(n-b-\frac{1}{2}\right)},
    \end{align*}
  we obtain \eqref{thm3.4-1} from \eqref{thm3.4-2}.     \qed

    Taking $b=-a, p=2m$ with $m\in\mathbb{N}$ in \eqref{thm3.4-1}, the following corollary holds.
    \begin{cor}\label{cor3}
    	For $m\in\mathbb{N}$ and $a\in\mathbb{C}$ with $a,a+\frac{1}{2}\notin\mathbb{Z}$, we have
    	\begin{align}\label{cor33}
    		&\sum_{n=1}^{\infty}\frac{(-1)^nh_n^{(2m)}}{\left(n+a-\frac{1}{2}\right)\left(n-a-\frac{1}{2}\right)}
            =\sum_{n=1}^{\infty}\frac{(-1)^nh_n^{(2m)}}{\left(n-\frac1{2}\right)^2-a^2} \notag\\
    		=&\frac{1}{2}\sum_{n=1}^{\infty}\frac{(-1)^n}{\left(n-\frac{1}{2}\right)^{2m}\left(\left(n-\frac1{2}\right)^2-a^2\right)}\notag\\
    		&+\frac{1}{2a}\sum_{k=0}^{m-1}\widetilde{t}(\overline{2k+1})\left(\overline{\zeta}(2m-2k; a)-\overline{\zeta}(2m-2k; -a)\right)\notag\\
    		&+\frac{\pi}{4a\cos(\pi a)}\left(\zeta(2m; a)-\zeta(2m; -a)\right).
    	\end{align}
    \end{cor}
    Setting $b=1-a, p=2m+1$ with $m\in\mathbb{N}_0$ in \eqref{thm3.4-1}, we get the following corollary.
    \begin{cor}\label{cor4}
    	For $m\in\mathbb{N}_0$ and $a\in\mathbb{C}$ with $a,a+\frac{1}{2}\notin\mathbb{Z}$, the following equality holds:
    	\begin{align}\label{cor44}
    		&\sum_{n=1}^{\infty}\frac{(-1)^nh_n^{(2m+1)}}{\left(n+a-\frac{1}{2}\right)\left(n-a+\frac{1}{2}\right)}
            =\sum_{n=1}^{\infty}\frac{(-1)^nh_n^{(2m+1)}}{n^2-\left(a-\frac1{2}\right)^2}\notag\\
    		=&\frac{1}{2a-1}\sum_{k=0}^m\widetilde{t}(\overline{2k+1})\left(\overline{\zeta}(2m-2k+1; 1-a)-\overline{\zeta}(2m-2k+1; a)\right)\notag\\
    		&+\frac{\pi}{2(2a-1)\cos(\pi a)}\left(2\widetilde{t}(2m+1)-\zeta(2m+1; a)-\zeta(2m+1; 1-a)\right).
    	\end{align}
    \end{cor}
	
More generally, we can get the following theorems.

\begin{thm}\label{thm-para1-ration-funct-residue}
Let $r(z)$ be a rational function which is $O(z^{-2})$ at infinity. Denote by $S$ the set of poles of $r(z)$. Assume that $0\notin S$ and $n,\pm(n-1/2)\not\in S$ for any $n\in\N$. Define
\[f(z):=\frac{\pi \tan(\pi z) \boldsymbol{\Psi}^{(p-1)}(1/2-z)}{(p-1)!}r(z).\]
Then for any $p\in \N$, we have
\begin{align}\label{para1-ration-funct-residue}
&-\sum_{n=1}^\infty \Big((-1)^p\t(p)+h_n^{(p)}\Big)r(n-1/2)-(-1)^p\sum_{n=1}^\infty \Big(\t(p)-h_{n-1}^{(p)}\Big)r(1/2-n)\nonumber\\
&+2\sum_{k=1}^{[p/2]} \frac{\t(2k)}{(p-2k)!}  \sum_{n=0}^\infty r^{(p-2k)}(n)+\sum_{\beta\in S}{\rm Res}(f(z),\beta)=0,
\end{align}
    where $\t(1)$ should be interpreted as $0$ wherever it occurs and $r^{(p)}(z)$ is the $p$-st derivative of $r(z)$.
\end{thm}	
	
\begin{thm}\label{thm-para1-ration-funct-residue-2}
Let $r(z)$ be a rational function which is $O(z^{-2})$ at infinity. Denote by $T$ the set of poles of $r(z)$. Assume that $0\notin T$ and $n,\pm(n-1/2)\not\in T$ for any $n\in\N$. Define
\[g(z):=\frac{\pi \boldsymbol{\Psi}^{(p-1)}(1/2-z)}{\cos(\pi z)(p-1)!}r(z).\]
Then for any $p\in \N$, we have
\begin{align}\label{para1-ration-funct-residue-2}
&\sum_{n=1}^\infty \Big((-1)^p\t(p)+h_n^{(p)}\Big)(-1)^n r(n-1/2)-(-1)^p\sum_{n=1}^\infty \Big(\t(p)-h_{n-1}^{(p)}\Big)(-1)^nr(1/2-n)\nonumber\\
&-2\sum_{k=0}^{[(p-1)/2]} \frac{\t(2k+1)}{(p-1-2k)!}  \sum_{n=0}^\infty (-1)^nr^{(p-1-2k)}(n)+\sum_{\beta\in T}{\rm Res}(g(z),\beta)=0,
\end{align}
where $\t(1)$ should be interpreted as $0$ wherever it occurs and $r^{(p)}(z)$ is the $p$-st derivative of $r(z)$.
\end{thm}		

    The proofs of Theorems \ref{thm-para1-ration-funct-residue} and \ref{thm-para1-ration-funct-residue-2} are completely similar as those of Theorems \ref{thm3.1} and \ref{thm3.4}. We leave the details to the interested readers.
    It is clear that if we set $r(z)=1/((z+a)(z+b))$, Theorems \ref{thm3.1} and \ref{thm3.4} follow immediately from Theorems \ref{thm-para1-ration-funct-residue} and \ref{thm-para1-ration-funct-residue-2}, respectively.

\begin{rem}
	For the (alternating) parametric nonlinear Euler $T$-sums, we can consider the following contour integral
	$$\oint\limits_{(\infty)} \frac{\pi \tan(\pi z) \boldsymbol{\Psi}^{(p_1-1)}\left(\frac{1}{2}-z\right)\boldsymbol{\Psi}^{(p_2-1)}\left(\frac{1}{2}-z\right)\cdots\boldsymbol{\Psi}^{(p_m-1)}\left(\frac{1}{2}-z\right)}{(p_1-1)!(p_2-1)!\cdots(p_m-1)!}r(z) \mathrm{d}z=0$$
	and
	$$\oint\limits_{(\infty)} \frac{\pi \boldsymbol{\Psi}^{(p_1-1)}\left(\frac{1}{2}-z\right)\boldsymbol{\Psi}^{(p_2-1)}\left(\frac{1}{2}-z\right)\cdots\boldsymbol{\Psi}^{(p_m-1)}\left(\frac{1}{2}-z\right)}{\cos(\pi z)(p_1-1)!(p_2-1)!\cdots(p_m-1)!}r(z) \mathrm{d}z=0.$$
\end{rem}


	\section{Double $t$-values and double $T$-values}\label{Sec:MtV}
	
	In \cite{BBB}, Borweins and Bradley used the evaluation of a specific parametirc Euler sum to obtain an explicit formula of double zeta values $\zeta(2j, 2m+1)$ by using power series expansion and comparing coeffcients. In this section, by using the same method, we give some explicit formulae for double $t$-values, double $T$-values, alternating double $t$-values and alternating double $T$-values.
	
	\begin{thm}\label{t-values-even-odd}
		For $j\in\mathbb{N}$ and $m\in\mathbb{N}_0$, we have
		\begin{align}\label{t-values-even-odd1}
			t(2j, 2m+1)=&t(2j)t(2m+1)-\frac{1}{2}t(2j+2m+1)\notag\\
			&-\sum_{k=1}^{m}\dbinom{2j+2m-2k}{2j-1}\frac{\zeta(2j+2m-2k+1)}{2^{2j+2m-2k+1}}t(2k)\notag\\
			&-\sum_{l=1}^{j}\dbinom{2j+2m-2l}{2m}\frac{\zeta(2j+2m-2l+1)}{2^{2j+2m-2l+1}}t(2l).
		\end{align}
	\end{thm}
	\proof
	By calculations, for $|a|<\frac{1}{2}$, we can easily deduce the following expansions
    \begin{align*}
	&\begin{aligned}
		\sum_{n=1}^{\infty}\frac{h_{n-1}^{(2m+1)}}{\left(n+a-\frac{1}{2}\right)\left(n-a-\frac{1}{2}\right)}=\sum_{i=1}^{\infty}\widetilde{t}(2i, 2m+1)a^{2i-2},
	\end{aligned}\\
	&\begin{aligned}
		\sum_{n=1}^{\infty}\frac{1}{\left(n-\frac{1}{2}\right)^{2m+1}\left(\left(n-\frac1{2}\right)^2-a^2\right)}=\sum_{i=1}^{\infty}\widetilde{t}(2i+2m+1)a^{2i-2},
	\end{aligned}\\
	&\begin{aligned}
		&\zeta(2m-2k+2; a)-\zeta(2m-2k+2; -a)\\
		=&-2\sum_{i=1}^{\infty}\dbinom{2i+2m-2k}{2i-1}\zeta(2i+2m+1-2k)a^{2i-1},
	\end{aligned}\\
    &\begin{aligned}
    	\zeta(2m+1; a)+\zeta(2m+1; -a)=2\sum_{i=0}^{\infty}\dbinom{2i+2m}{2i}\zeta(2i+2m+1)a^{2i}.
    \end{aligned}
    	\end{align*}
	Substituting these identities and \eqref{tan} into \eqref{cor11}, we get
	\begin{align*}
		\sum_{i=1}^{\infty}\widetilde{t}(2i, 2m+1)a^{2i-2}=&-\frac{1}{2}\sum_{i=1}^{\infty}\widetilde{t}(2i+2m+1)a^{2i-2}\\
		&-\sum_{i=1}^{\infty}\sum_{k=1}^m\dbinom{2i+2m-2k}{2i-1}\zeta(2i+2m+1-2k)\widetilde{t}(2k)a^{2i-2}\\
		&+\sum_{i=1}^{\infty}\widetilde{t}(2m+1)\widetilde{t}(2i)a^{2i-2}\\
		&-\sum_{i=1}^{\infty}\sum_{l=1}^{i}\dbinom{2i+2m-2l}{2m}\zeta(2i+2m+1-2l)\widetilde{t}(2l)a^{2i-2}.
	\end{align*}
    By comparing the coefficients of $a^{2j-2}$ in the above equation and using \eqref{widetilde-t}, we obtain \eqref{t-values-even-odd1}. This completes the proof.
    \qed

    \begin{thm}\label{t-values-odd-even}
    	For $j,m\in\mathbb{N}$, the following identity holds:
    	\begin{align}\label{t-values-odd-even1}
    		t(2j+1, 2m)=&-\frac{1}{2}t(2j+2m+1)\notag\\
    		&+\sum_{k=1}^{m}\dbinom{2j+2m-2k}{2j}\frac{\zeta(2j+2m-2k+1)}{2^{2j+2m-2k+1}}t(2k)\notag\\
    		&+\sum_{l=1}^{j}\dbinom{2j+2m-2l}{2m-1}\frac{\zeta(2j+2m-2l+1)}{2^{2j+2m-2l+1}}t(2l).
    	\end{align}
    \end{thm}
	\proof
	By taking the first derivative of both sides of \eqref{thm3.1} with respect to $a$ and setting $b=-a, p=2m$, we deduce that
	\begin{align}\label{t-values-odd-even2}
		&\sum_{n=1}^{\infty}\frac{h_n^{(2m)}}{\left(n+a-\frac{1}{2}\right)^2\left(n-a-\frac{1}{2}\right)}+
		\sum_{n=1}^{\infty}\frac{h_{n-1}^{(2m)}}{\left(n+a-\frac{1}{2}\right)\left(n-a-\frac{1}{2}\right)^2}\notag\\
		=&-\frac{1}{a}\sum_{k=1}^m\widetilde{t}(2k)(2m-2k+1)\zeta(2m-2k+2; a)\notag\\
		&-\frac{1}{2a^2}\sum_{k=1}^m\widetilde{t}(2k)\left(\zeta(2m-2k+1; a)-\zeta(2m-2k+1; -a)\right)\notag\\
		&+\frac{1}{4a^2}\pi \tan(\pi a)\left(\zeta(2m; a)+\zeta(2m; -a)-2\widetilde{t}(2m)\right)\notag\\
		&+\frac{1}{2a}\left\{2m\pi \tan(\pi a)\zeta(2m+1;a)-(\zeta(2m;a)-\widetilde{t}(2m))\frac{d}{da}(\pi \tan(\pi a))\right\}.
	\end{align}
	The following process is similar to the proof of Theorem \ref{t-values-even-odd}. By direct calculations, for $|a|<\frac{1}{2}$, we have
    \begin{align*}
	&\begin{aligned}
	    \sum_{n=1}^{\infty}\frac{h_{n-1}^{(2m)}}{\left(n+a-\frac{1}{2}\right)^2\left(n-a-\frac{1}{2}\right)}
	    &=\sum_{i=1}^{\infty}\left(\sum_{l=1}^{i-1}(-1)^{l-1}l\right)\widetilde{t}(i+1, 2m)a^{i-2}\\
	    &=\sum_{i=1}^{\infty}(-1)^i\left[\frac{i}{2}\right]\widetilde{t}(i+1, 2m)a^{i-2},
	\end{aligned}\\
	&\begin{aligned}
		\sum_{n=1}^{\infty}\frac{h_{n-1}^{(2m)}}{\left(n+a-\frac{1}{2}\right)\left(n-a-\frac{1}{2}\right)^2}
		=\sum_{i=1}^{\infty}\left[\frac{i}{2}\right]\widetilde{t}(i+1, 2m)a^{i-2},
	\end{aligned}\\
    &\begin{aligned}
        \sum_{n=1}^{\infty}\frac{1}{\left(n-\frac{1}{2}\right)^{2m}\left(n+a-\frac{1}{2}\right)^2\left(n-a-\frac{1}{2}\right)}=\sum_{i=1}^{\infty}(-1)^i\left[\frac{i}{2}\right]\widetilde{t}(i+2m+1)a^{i-2},
    \end{aligned}\\
    &\begin{aligned}
    	&\zeta(2m-2k+1; a)-\zeta(2m-2k+1; -a)\\
    	=&\frac{2}{a^{2m-2k+1}}-2\sum_{i=1}^{\infty}\dbinom{2i+2m-2k-1}{2m-2k}\zeta(2i+2m-2k)a^{2i-1},
    \end{aligned}\\
    &\begin{aligned}
    	\zeta(2m; a)+\zeta(2m; -a)=\frac{2}{a^{2m}}+2\sum_{i=1}^{\infty}\dbinom{2i+2m-3}{2i-2}\zeta(2i+2m-2)a^{2i-2}.
    \end{aligned}
     \end{align*}
	The power series expansion of the Hurwitz zeta function is
	\begin{align}\label{hurwitz-expansion}
		\zeta(z;a)=\frac{1}{a^{z}}+\sum_{i=1}^{\infty}(-1)^{i-1}\dbinom{i+z-2}{i-1}\zeta(i+z-1)a^{i-1}.
	\end{align}
	Then by substituting the above formulae and \eqref{tan} into \eqref{t-values-odd-even2}, we arrive at
	\begin{align}\label{t-value-odd-even3}
		&\sum_{i=1}^{\infty}(1+(-1)^i)\left[\frac{i}{2}\right]\widetilde{t}(i+1, 2m)a^{i-2}\notag\\
		=&\sum_{i=1}^{\infty}(-1)^{i-1}\left[\frac{i}{2}\right]\widetilde{t}(i+2m+1)a^{i-2}-\sum_{k=1}^{m}(2m-2k+2)\widetilde{t}(2k)a^{2k-2m-3}\notag\\
		&+\sum_{i=1}^{\infty}(-1)^i\sum_{k=1}^{m}\dbinom{i+2m-2k}{i-1}(2m-2k+1)\zeta(i+2m-2k+1)\widetilde{t}(2k)a^{i-2}\notag\\
		&+\sum_{i=1}^{\infty}\sum_{k=1}^{m}\dbinom{2i+2m-2k-1}{2m-2k}\zeta(2i+2m-2k)\widetilde{t}(2k)a^{2i-3}\notag\\
		&+\sum_{i=1}^{\infty}(2i-2)\widetilde{t}(2m)\widetilde{t}(2i)a^{2i-3}+\sum_{i=1}^{\infty}(2m-2i+2)\widetilde{t}(2i)a^{2i-2m-3}\notag\\
		&+\sum_{i=1}^{\infty}\sum_{l=1}^{i-1}\dbinom{2i+2m-2l-3}{2i-2l-1}\zeta(2i+2m-2l-2)\widetilde{t}(2l)a^{2i-5}\notag\\
		&+\sum_{i=1}^{\infty}\sum_{l=1}^{\left[\frac{i-1}{2}\right]}(-1)^{i-1}\dbinom{i+2m-2l-1}{2m}(2m)\zeta(i+2m-2l)\widetilde{t}(2l)a^{i-3}\notag\\
		&+\sum_{i=1}^{\infty}\sum_{l=1}^{\left[\frac{i-1}{2}\right]}(-1)^{i}\dbinom{i+2m-2l-2}{2m-1}(2l-1)\zeta(i+2m-2l-1)\widetilde{t}(2l)a^{i-4}.
	\end{align}
	By comparing the coefficients of $a^{2j-2}$ in \eqref{t-value-odd-even3} and applying \eqref{widetilde-t}, we deduce \eqref{t-values-odd-even1}.
	\qed
	
	In the same manner we prove the following formulae for the alternating double $t$-values.
		\begin{thm}\label{alternating-t-values-even}
		For $j, m\in\mathbb{N}$, we have
		\begin{align}\label{alternating-t-values-even1}
			t(\overline{2j}, 2m)=&-\frac{1}{2}t(\overline{2j+2m})\notag\\
			&+\sum_{k=0}^{m-1}\dbinom{2j+2m-2k-2}{2j-1}\frac{\overline{\zeta}(2j+2m-2k-1)}{2^{2j+2m-2k-1}}t(\overline{2k+1})\notag\\
			&+\sum_{l=0}^{j-1}\dbinom{2j+2m-2l-2}{2m-1}\frac{\zeta(2j+2m-2l-1)}{2^{2j+2m-2l-1}}t(\overline{2l+1}).
		\end{align}
	\end{thm}
	\proof
	By calculations, for $|a|<\frac{1}{2}$, we deduce the following expansions
    \begin{align*}
	&\begin{aligned}
		\sum_{n=1}^{\infty}\frac{(-1)^nh_{n-1}^{(2m)}}{\left(n+a-\frac{1}{2}\right)\left(n-a-\frac{1}{2}\right)}=\sum_{i=1}^{\infty}\widetilde{t}(\overline{2i}, 2m)a^{2i-2},
	\end{aligned}\\
	&\begin{aligned}
	\sum_{n=1}^{\infty}\frac{(-1)^n}{\left(n-\frac{1}{2}\right)^{2m}\left(\left(n-\frac1{2}\right)^2-a^2\right)}=\sum_{i=1}^{\infty}\widetilde{t}(\overline{2i+2m})a^{2i-2},
	\end{aligned}\\
	&\begin{aligned}
		&\overline{\zeta}(2m-2k; a)-\overline{\zeta}(2m-2k; -a)\\
		=&2\sum_{i=1}^{\infty}\dbinom{2i+2m-2k-2}{2i-1}\overline{\zeta}(2i+2m-2k-1)a^{2i-1},
	\end{aligned}\\
	&\begin{aligned}
		\zeta(2m; a)-\zeta(2m; -a)=-2\sum_{i=1}^{\infty}\dbinom{2i+2m-2}{2i-1}\zeta(2i+2m-1)a^{2i-1}.
	\end{aligned}
    \end{align*}
	Substituting these identities and \eqref{cos} into \eqref{cor33}, we obtain
	\begin{align*}
		\sum_{i=1}^{\infty}\widetilde{t}(\overline{2i}, 2m)a^{2i-2}=
		&-\frac{1}{2}\sum_{i=1}^{\infty}\widetilde{t}(\overline{2i+2m})a^{2i-2}\\
		&+\sum_{i=1}^{\infty}\sum_{k=0}^{m-1}\dbinom{2i+2m-2k-2}{2i-1}\overline{\zeta}(2i+2m-2k-1)\widetilde{t}(\overline{2k+1})a^{2i-2}\\
		&+\sum_{i=1}^{\infty}\sum_{l=0}^{i-1}\dbinom{2i+2m-2l-2}{2m-1}\zeta(2i+2m-2l-1)\widetilde{t}(\overline{2l+1})a^{2i-2}.
	\end{align*}
	By comparing the coefficients of $a^{2j-2}$ in the above equation and using \eqref{widetilde-t}, we derive the desired result. This completes the proof of Theorem \ref{alternating-t-values-even}.
	\qed
	
	\begin{thm}\label{alternating-t-values-odd}
		For $j, m\in\mathbb{N}_0$, the following identity holds:
		\begin{align}\label{alternating-t-values-odd1}
			t(\overline{2j+1}, 2m+1)=&t(\overline{2j+1})t(2m+1)-\frac{1}{2}t(\overline{2j+2m+2})\notag\\
			&-\sum_{k=0}^{m}\dbinom{2j+2m-2k}{2j}\frac{\overline{\zeta}(2j+2m-2k+1)}{2^{2j+2m-2k+1}}t(\overline{2k+1})\notag\\
			&-\sum_{l=0}^{j}\dbinom{2j+2m-2l}{2m}\frac{\zeta(2j+2m-2l+1)}{2^{2j+2m-2l+1}}t(\overline{2l+1}).
		\end{align}
	\end{thm}
	\proof
	By taking the first derivative of both sides of \eqref{thm3.4-1} with respect to $a$ and setting $b=-a, p=2m+1$, we conclude that
	\begin{align}\label{alternating-t-values-odd2}
		&\sum_{n=1}^{\infty}\frac{(-1)^nh_n^{(2m+1)}}{\left(n+a-\frac{1}{2}\right)^2\left(n-a-\frac{1}{2}\right)}+
		\sum_{n=1}^{\infty}\frac{(-1)^nh_{n-1}^{(2m+1)}}{\left(n+a-\frac{1}{2}\right)\left(n-a-\frac{1}{2}\right)^2}\notag\\
		=&-\frac{1}{a}\sum_{k=0}^m\widetilde{t}(\overline{2k+1})(2m-2k+1)\overline{\zeta}(2m-2k+2; a)\notag\\
		&-\frac{1}{2a^2}\sum_{k=0}^m\widetilde{t}(\overline{2k+1})\left(\overline{\zeta}(2m-2k+1; a)-\overline{\zeta}(2m-2k+1; -a)\right)\notag\\
		&+\frac{1}{4a^2}\frac{\pi}{\cos(\pi a)}\left(\zeta(2m+1; -a)-\zeta(2m+1; a)\right)\notag\\
		&-\frac{1}{2a}\left\{\frac{\pi}{\cos(\pi a)}(2m+1)\zeta(2m+2;a)-(\zeta(2m+1;a)-\widetilde{t}(2m+1))\frac{d}{da}\left(\frac{\pi}{\cos(\pi a)}\right)\right\}.
	\end{align}
    By direct calculations, for $|a|<\frac{1}{2}$, we have
    \begin{align*}
	&\begin{aligned}
		\sum_{n=1}^{\infty}\frac{(-1)^nh_{n-1}^{(2m+1)}}{\left(n+a-\frac{1}{2}\right)^2\left(n-a-\frac{1}{2}\right)}
		=\sum_{i=1}^{\infty}(-1)^i\left[\frac{i}{2}\right]\widetilde{t}(\overline{i+1}, 2m+1)a^{i-2},
	\end{aligned}\\
	&\begin{aligned}
		\sum_{n=1}^{\infty}\frac{(-1)^nh_{n-1}^{(2m+1)}}{\left(n+a-\frac{1}{2}\right)\left(n-a-\frac{1}{2}\right)^2}
		=\sum_{i=1}^{\infty}\left[\frac{i}{2}\right]\widetilde{t}(\overline{i+1}, 2m+1)a^{i-2},
	\end{aligned}\\
	&\begin{aligned} \sum_{n=1}^{\infty}\frac{(-1)^n}{\left(n-\frac{1}{2}\right)^{2m+1}\left(n+a-\frac{1}{2}\right)^2\left(n-a-\frac{1}{2}\right)}=\sum_{i=1}^{\infty}(-1)^i\left[\frac{i}{2}\right]\widetilde{t}(\overline{i+2m+2})a^{i-2},
	\end{aligned}\\
	&\begin{aligned}
		&\overline{\zeta}(2m-2k+1; a)-\overline{\zeta}(2m-2k+1; -a)\\
		=&\frac{2}{a^{2m-2k+1}}+2\sum_{i=1}^{\infty}\dbinom{2i+2m-2k-1}{2m-2k}\overline{\zeta}(2i+2m-2k)a^{2i-1},
	\end{aligned}\\
	&\begin{aligned}
		\zeta(2m+1; -a)-\zeta(2m+1; a)=-\frac{2}{a^{2m+1}}+2\sum_{i=1}^{\infty}\dbinom{2i+2m-1}{2m}\zeta(2i+2m)a^{2i-1},
	\end{aligned}\\
	&\begin{aligned}
	    \overline{\zeta}(2m-2k+2;a)=\frac{1}{a^{2m-2k+2}}+\sum_{i=1}^{\infty}(-1)^{i}\dbinom{i+2m-2k}{i-1}\overline{\zeta}(i+2m-2k+1)a^{i-1}.
    \end{aligned}
       \end{align*}
	Then by substituting the above formulae and Eqs. \eqref{cos}, \eqref{hurwitz-expansion} into \eqref{alternating-t-values-odd2}, we deduce that
	\begin{align}\label{alternating-t-values-odd3}
		&\sum_{i=1}^{\infty}(1+(-1)^i)\left[\frac{i}{2}\right]\widetilde{t}(\overline{i+1}, 2m+1)a^{i-2}\notag\\
		=&\sum_{i=1}^{\infty}(-1)^{i-1}\left[\frac{i}{2}\right]\widetilde{t}(\overline{i+2m+2})a^{i-2}-\sum_{k=0}^{m}(2m-2k+2)\widetilde{t}(\overline{2k+1})a^{2k-2m-3}\notag\\
		&-\sum_{i=1}^{\infty}(-1)^i\sum_{k=0}^{m}\dbinom{i+2m-2k}{i-1}(2m-2k+1)\overline{\zeta}(i+2m-2k+1)\widetilde{t}(\overline{2k+1})a^{i-2}\notag\\
		&-\sum_{i=1}^{\infty}\sum_{k=0}^{m}\dbinom{2i+2m-2k-1}{2m-2k}\overline{\zeta}(2i+2m-2k)\widetilde{t}(\overline{2k+1})a^{2i-3}\notag\\
		&+\sum_{i=0}^{\infty}(2i)\widetilde{t}(2m+1)\widetilde{t}(\overline{2i+1})a^{2i-2}+\sum_{i=0}^{\infty}(2m-2i+2)\widetilde{t}(\overline{2i+1})a^{2i-2m-3}\notag\\
		&-\sum_{i=1}^{\infty}\sum_{l=0}^{i-1}\dbinom{2i+2m-2l-1}{2m}\zeta(2i+2m-2l)\widetilde{t}(\overline{2l+1})a^{2i-3}\notag\\
		&-\sum_{i=1}^{\infty}\sum_{l=0}^{\left[\frac{i-1}{2}\right]}(-1)^{i}\dbinom{i+2m-2l}{2m+1}(2m+1)\zeta(i+2m-2l+1)\widetilde{t}(\overline{2l+1})a^{i-2}\notag\\
		&+\sum_{i=1}^{\infty}\sum_{l=0}^{\left[\frac{i-1}{2}\right]}(-1)^{i}\dbinom{i+2m-2l-1}{2m}(2l)\zeta(i+2m-2l)\widetilde{t}(\overline{2l+1})a^{i-3}.
	\end{align}
	Thus, comparing the coefficients of $a^{2j-2}$ in \eqref{alternating-t-values-odd3} and applying \eqref{widetilde-t}, we derive the desired result.
	\qed

	Proceeding in a similar manner to the proofs of Theorem \ref{t-values-even-odd}-\ref{alternating-t-values-odd}, we obtain the explicit formulae for double $T$-values and alternating double $T$-values.
	
	\begin{thm}\label{T-values-even-odd}
		For $j\in\mathbb{N}$, $m\in\mathbb{N}_0$, we have
		\begin{align}\label{T-values-even-odd1}
			T(2j,2m+1)=&\dbinom{2m+2j}{2m}T(2m+2j+1)\notag\\
			&-\sum_{k=1}^{m}\dbinom{2m+2j-2k}{2j-1}T(2m+2j-2k+1)T(2k)\notag\\
			&-\sum_{l=1}^{j-1}\dbinom{2m+2j-2l}{2m}\frac{\zeta(2l)}{2^{2l-1}}T(2m+2j-2l+1).
		\end{align}
	\end{thm}
	\proof
	Replacing $a$ by $a+\frac{1}{2}$ and letting $b=-a+\frac{1}{2}$ and $p=2m+1$ in \eqref{thm3.1-1}, we deduce that
	\begin{align}\label{T-values-even-odd2}
		&\sum_{n=1}^{\infty}\frac{h_n^{(2m+1)}}{(n+a)(n-a)}\notag\\
		=&\frac{1}{2a}\sum_{k=1}^{m}\widetilde{t}(2k)\left\{\zeta\left(2m-2k+2;a+\frac{1}{2}\right)-\zeta\left(2m-2k+2;-a+\frac{1}{2}\right)\right\}\notag\\
		&+\frac{1}{4a}\pi \cot(\pi a)\left\{\zeta\left(2m+1;a+\frac{1}{2}\right)+\zeta\left(2m+1;-a+\frac{1}{2}\right)-2\widetilde{t}(2m+1)\right\}.
	\end{align}
	According to definitions and direct calculations, we have
	\begin{align*}
	&\begin{aligned}
    	\pi \cot(\pi a)=\frac{1}{a}-2\sum_{i=1}^{\infty}\zeta(2i)a^{2i-1},
	\end{aligned}\\
	&\begin{aligned}
	   	T(s_1,s_2)=\frac{1}{2^{s_1+s_2-2}}\sum_{n=1}^{\infty}\frac{h_n^{(s_2)}}{n^{s_1}},
    \end{aligned}\\
    &\begin{aligned}
    	\widetilde{t}(s)=2^{s-1}T(s).
    \end{aligned}
	\end{align*}
	By applying these equations and calculating the power series expansion of \eqref{T-values-even-odd2}, we can obtain \eqref{T-values-even-odd1} by comparing the coefficients of $a^{2j-2}$. This completes the proof.
	\qed
	
	\begin{thm}\label{T-values-odd-even}
		For $j, m\in\mathbb{N}$ the following identity holds:
		\begin{align}\label{T-values-odd-even1}
			T(2j+1,2m)=&-\dbinom{2m+2j}{2j+1}T(2m+2j+1)\notag\\
			&+\sum_{k=1}^{m}\dbinom{2m+2j-2k}{2j}T(2m+2j-2k+1)T(2k)\notag\\
			&+\sum_{l=1}^{j}\dbinom{2m+2j-2l}{2m-1}\frac{\zeta(2l)}{2^{2l-1}}T(2m+2j-2l+1).
		\end{align}
	\end{thm}
	\proof
	By taking the first derivative of both sides of \eqref{thm3.1} with respect to $a$ and replacing $a$ by $a+\frac{1}{2}$ and setting $b=-a+\frac{1}{2}, p=2m$, we arrive at
	\begin{align}\label{T-values-odd-even2}
		&\sum_{n=1}^{\infty}\frac{h_n^{(2m)}}{(n+a)^2(n-a)}+\sum_{n=1}^{\infty}\frac{h_n^{(2m)}}{(n+a)(n-a)^2}\notag\\
		=&\frac{1}{a}\sum_{k=1}^{m}\widetilde{t}(2k)(2k-2m-1)\zeta\left(2m-2k+2;a+\frac{1}{2}\right)\notag\\
		&-\frac{1}{2a^2}\sum_{k=1}^{m}\widetilde{t}(2k)\left\{\zeta\left(2m-2k+1;a+\frac{1}{2}\right)-\zeta\left(2m-2k+1;-a+\frac{1}{2}\right)\right\}\notag\\
		&-\frac{1}{4a^2}\pi \cot(\pi a)\left\{\zeta\left(2m;a+\frac{1}{2}\right)+\zeta\left(2m;-a+\frac{1}{2}\right)-2\widetilde{t}(2m)\right\}\notag\\
		&+\frac{1}{2a}\left\{\left(\zeta\left(2m;a+\frac{1}{2}\right)-\widetilde{t}(2m)\right)\frac{d}{da}(\pi \cot(\pi a))-2m\pi \cot(\pi a)\zeta\left(2m+1;a+\frac{1}{2}\right)\right\}.
	\end{align}
	Thus, calculating the power series expansion of \eqref{T-values-odd-even2} and comparing the coefficients of $a^{2j-2}$, we derive the desired result.
	\qed
	
	\begin{thm}\label{alternating-T-values-even-odd}
		For $j\in\mathbb{N}$ and $m\in\mathbb{N}_0$, we have
		\begin{align}\label{alternating-T-values-even-odd1}
			T(\overline{2j},2m+1)=&-\dbinom{2m+2j}{2m}T(2m+2j+1)\notag\\
			&+\sum_{k=0}^{m}\dbinom{2m+2j-2k-1}{2j-1}T(\overline{2m+2j-2k})T(\overline{2k+1})\notag\\
			&-\sum_{l=1}^{j-1}\dbinom{2m+2j-2l}{2m}\frac{\overline{\zeta}(2l)}{2^{2l-1}}T(2m+2j-2l+1).
		\end{align}
	\end{thm}
	\proof
		Replacing $a$ by $a+\frac{1}{2}$ and letting $b=-a+\frac{1}{2}$ and $p=2m+1$ in \eqref{thm3.4-1}, we deduce that
	\begin{align}\label{alternating-T-values-even-odd2}
		&\sum_{n=1}^{\infty}\frac{(-1)^nh_n^{(2m+1)}}{(n+a)(n-a)}\notag\\
		=&\frac{1}{2a}\sum_{k=0}^{m}\widetilde{t}(\overline{2k+1})\left\{\overline{\zeta}\left(2m-2k+1;-a+\frac{1}{2}\right)-\overline{\zeta}\left(2m-2k+1;a+\frac{1}{2}\right)\right\}\notag\\
		&+\frac{1}{4a}\frac{\pi}{\sin(\pi a)}\left\{\zeta\left(2m+1;a+\frac{1}{2}\right)+\zeta\left(2m+1;-a+\frac{1}{2}\right)-2\widetilde{t}(2m+1)\right\}.
	\end{align}
	According to definitions and direct calculations, we find that
	\begin{align*}
		&\begin{aligned}
			\frac{\pi}{\sin(\pi a)}=\frac{1}{a}+2\sum_{i=1}^{\infty}\overline{\zeta}(2i)a^{2i-1},
		\end{aligned}\\
		&\begin{aligned}
			T(\overline{s_1},s_2)=\frac{1}{2^{s_1+s_2-2}}\sum_{n=1}^{\infty}\frac{(-1)^{n-1}h_n^{(s_2)}}{n^{s_1}},
		\end{aligned}\\
		&\begin{aligned}
			\widetilde{t}(\overline{s})=2^{s-1}T(\overline{s}).
		\end{aligned}
	\end{align*}
	Applying these equations and calculating the power series expansion of \eqref{alternating-T-values-even-odd2}, we obtain \eqref{alternating-T-values-even-odd1} by comparing the coefficients of $a^{2j-2}$.
	\qed
	
	\begin{thm}\label{alternating-T-values-odd-even}
		For $j\in\mathbb{N}_0$, $m\in\mathbb{N}$, the following identity holds:
		\begin{align}\label{alternating-T-values-odd-even1}
			T(\overline{2j+1},2m)=&\dbinom{2m+2j}{2j+1}T(2m+2j+1)\notag\\
			&-\sum_{k=0}^{m-1}\dbinom{2m+2j-2k-1}{2j}T(\overline{2m+2j-2k})T(\overline{2k+1})\notag\\
			&+\sum_{l=1}^{j}\dbinom{2m+2j-2l}{2m-1}\frac{\overline{\zeta}(2l)}{2^{2l-1}}T(2m+2j-2l+1).
		\end{align}
	\end{thm}
	\proof
	By taking the first derivative of both sides of \eqref{thm3.4-1} with respect to $a$ and replacing $a$ by $a+\frac{1}{2} $ and setting $b=-a+\frac{1}{2}, p=2m$, we arrive at
	\begin{align}\label{alternating-T-values-odd-even2}
		&\sum_{n=1}^{\infty}\frac{(-1)^{n-1}h_n^{(2m)}}{(n+a)^2(n-a)}+\sum_{n=1}^{\infty}\frac{(-1)^{n-1}h_n^{(2m)}}{(n+a)(n-a)^2}\notag\\
		=&-\frac{1}{a}\sum_{k=0}^{m-1}\widetilde{t}(\overline{2k+1})(2m-2k)\overline{\zeta}\left(2m-2k+1;a+\frac{1}{2}\right)\notag\\
		&+\frac{1}{2a^2}\sum_{k=0}^{m-1}\widetilde{t}(\overline{2k+1})\left\{\overline{\zeta}\left(2m-2k;-a+\frac{1}{2}\right)-\overline{\zeta}\left(2m-2k;a+\frac{1}{2}\right)\right\}\notag\\
		&+\frac{1}{4a^2}\frac{\pi}{\sin(\pi a)}\left\{\zeta\left(2m;a+\frac{1}{2}\right)+\zeta\left(2m;-a+\frac{1}{2}\right)-2\widetilde{t}(2m)\right\}\notag\\
		&-\frac{1}{2a}\left\{\left(\zeta\left(2m;a+\frac{1}{2}\right)-\widetilde{t}(2m)\right)\frac{d}{da}\left(\frac{\pi}{\sin(\pi a)}\right)-2m\frac{\pi}{\sin(\pi a)}\zeta\left(2m+1;a+\frac{1}{2}\right)\right\}.
	\end{align}
	By calculating the power series expansion of \eqref{alternating-T-values-odd-even2} and comparing the coefficients of $a^{2j-2}$, we obtain \eqref{alternating-T-values-odd-even1}. This completes the proof of Theorem \ref{alternating-T-values-odd-even}.
	\qed

\vskip10pt
	
\noindent{\bf Acknowledgement.} The authors express the deepest gratitude to their supervisor Professor
Zhonghua Li for his valuable comments and encouragement. Ce Xu is supported by the National Natural Science Foundation of China (Grant No. 12101008), the Natural Science Foundation of Anhui Province (Grant No. 2108085QA01) and the University Natural Science Research Project of Anhui Province (Grant No. KJ2020A0057). Lu Yan is supported by the Fundamental Research Funds for the Central Universities.

\medskip


\end{document}